\documentclass[10pt]{article}
\usepackage{subfigure}
\usepackage{graphicx}
\usepackage{amssymb}
\usepackage{amsmath}
\usepackage{float,color}
\usepackage{amsthm}
\usepackage{mathrsfs}
\newtheorem{theorem}{Theorem}[section]

\newtheorem{lemma}{Lemma}[section]


\numberwithin{equation}{section}

\def\NN{\mathbb{N}}

\begin{document}
\title{\bf Performance of Lagrangian Descriptors and Their Variants in Incompressible Flows}
\author{{Alfonso Ruiz-Herrera}\footnote{Departamento de
Matem\'{a}ticas, Universidad de Oviedo, Spain (alfonsoruiz@dma.uvigo.es, ruizalfonso@uniovi.es).}
}

\maketitle

\begin{abstract}
{The method of Lagrangian Descriptors (LDs) has been  applied in many different contexts, specially in geophysical flows.
In this paper we  analyze their performance   in incompressible flows.   We construct broad families of systems where  this diagnostic fails in the  detection of barriers to transport.  Another aim of this manuscript is to illustrate the same deficiencies  in the recent diagnostic proposed by Craven and Hern\'{a}ndez.}
\end{abstract}

\maketitle

{\bf The method of Lagrangian Descriptors (LDs) is a procedure to detect barriers to transport. Originally proposed by Mancho, Wiggins and their co-workers, this method has a basic heuristic law:  invariant manifolds (hyperbolic structures) are detected via singular features of the $M$-function.
This function is defined in terms of the arc length of the orbits of the system. In this paper we present a large family of incompressible flows in which
LDs do not detect these structures. Actually, we observe that in some cases all singular features are evanescent at the invariant manifolds. Although this method has been intensively applied in geophysical flows, the conclusion is that it is not founded on solid mathematical grounds. Generally speaking, the singularities associated with the $M$-function in a system are unrelated to the detection of invariant manifolds. Recently, Craven and Hern\'{a}ndez have  extended the method of  LDs to  the context of thermalized chemical reactions. In the second part of this manuscript we generalize the previous conclusions to this setting.}

\section{Introduction}

The method of Lagrangian Descriptors (LDs)
is a diagnostic to  detect barriers to transport in fluid dynamics or  invariant manifolds in dynamical systems \cite{mendoza-mancho PRL}-\cite{lopesino}. This method has been mainly applied by Mancho, Wiggins, and their co-workers  in different contexts involving transport phenomena \cite{mendoza-mancho PRL}-\cite{wiggins}. For instance,  Mendoza and Mancho  studied the skeleton of the flow for the Kuroshio current \cite{mendoza npg2012} or de la C\'{a}mara {\it et al.}  provided barriers to transport in the Antarctic polar vortex \cite{camara 1}.  The conclusions of those papers are nevertheless doubtful because there is no mathematical foundation of LDs beyond  heuristic arguments in simple systems. In this paper, we
 provide some families of incompressible flows where LDs do not detect the expected structures.

LDs  were originally designed to compute separatrices in flows with general dependence on time, a computation typically performed via Finite Time Lyapunov Exponents or the methods of Haller and his co-workers (variational and geodesic theories of Lagrangian Coherent Structures) \cite{haller}. The heuristic law of LDs is that invariant manifolds are associated with some singular features of the contour-lines of the $M$-function (this function is defined in terms of the arc length of the trajectories of the system).
Generally speaking, such law seems questionable. For instance, a prerequisite for reliable  predictions is the independence of the observer or objectivity \cite{chaos intro}. Unfortunately, as Haller nicely emphasized \cite{haller M, haller}, LDs are not objective.

In a previous work \cite{ruizherrera}, we  presented several pathologies and
 examples where LDs do not detect barriers to transport. Although those examples included Hamiltonian systems, they were rather involved and non-autonomous, (with an aperiodic time dependence). One of the aims of this paper is to present simple families of incompressible flows for which LDs do not detect  the invariant manifolds.  Note that the performance of LDs in higher dimension is rather limited. The method is already questionable in
\begin{equation}\label{S25I}\left\{\begin{array}{lll}
x_{1}'=0.5 x_{1}\\
x_{2}'= 1.5 x_{2}\\
x'_{3}=-2 x_{3}.\\
\end{array}\right.\end{equation}
In this system, the contour-lines of the $M$-function (for large  values of $\tau$) are hyperplanes parallel to $\{x_{3}=0\}$. Thus there is no established link between the dynamical skeleton of (\ref{S25I}) and the structure of the $M$-function over finite time.

 Recently, Lopesino {\it et al.}  \cite{lopesino}  have obtained some theoretical results in order to justify the performance  of LDs in
\begin{equation}\label{S2512}\left\{\begin{array}{lll}
x_{n+1}= \lambda x_{n}\\
y_{n+1}= \frac{1}{\lambda}y_{n}\;\;\;\;{\rm with}\;\;\lambda>1.
\end{array}\right.\end{equation}
 The analysis of  (\ref{S2512}) can be obtained with elementary tools but the application of LDs seems to suggest that the results in \cite{lopesino} could be extended to larger classes of systems. In the next section we conclude that this is not the case because the method does not work for a small perturbation of system (\ref{S2512}).

Another  variant of the method of LDs has been presented by Cranven and Hern\'{a}ndez to compute the invariant manifolds in thermalized chemical reactions \cite{craven and hernandez}-\cite{craven and hernandez 1}. This new approach, reminiscent to the original paper by Madrid and Mancho \cite{mancho chaos},  consists in minimizing the arc length of trajectories forwards in time to detect the stable manifold. The objective of this technique was to describe manifolds that separate reactive and non-reactive trajectories in the Langevin equation \cite{craven and hernandez}-\cite{craven and hernandez 1}. Again, the presentation is in a heuristic style and the theoretical aspects are not developed. We will describe some simple examples in the context of the Langevin equation where  invariant manifolds cannot be found via their method.

 An outline of the paper is as follows: Section 2 reviews the literature related to the method of LDs for the reader's convenience. Section 3 deals with the examples  of  incompressible flows. The analysis is restricted to this class of systems but our results can be applied to a broader range of dynamical systems, in particular, in most linear systems. Section 4 deals with the Langevin equation in connection with the method employed by Craven and Hern\'{a}ndez. Finally, we discuss the implications of our analysis in practical situations.

 \section{Lagrangian Descriptors: Definitions and Mathematical Foundation}

 LDs aim to describe the global dynamical picture of the geometrical structures for dynamical systems. In particular, they
 were used to detect stable and unstable manifolds and elliptic regions. This method in continuous systems typically involves the map defined as
 $$M(\tau; x_{0},t_{0})=\int_{t_{0}-\tau}^{t_{0}+\tau}\|v(x(t;t_{0},x_{0}),t)\|dt$$
where $x(t;t_{0},x_{0})$ is the unique solution satisfying $x(t_{0})=x_{0}$ for a smooth velocity field
\begin{equation}\label{eq}
x'(t)=v(x(t),t).
\end{equation}
The $M$-function was introduced  by Madrid and Mancho in \cite{mancho chaos} and  represents the arc length described by the trajectory starting at $x(t_{0})$ when it evolves forwards and backwards in time for a period $\tau$. This map is smooth  except at equilibria. In discrete planar systems, we define for each $N\in\NN$,
\begin{equation}\label{Mp}
 MD_{p}(x_{0},y_{0})=\sum_{i=-N}^{N-1}|x_{i+1}-x_{i}|^{p}+|y_{i+1}-y_{i}|^{p}\;\;\;\;\;{\rm with}\;0<p<1.
 \end{equation}
with $\{x_{n}, y_{n}\}_{n=-N}^{n=N}$   an orbit of length $2N+1$ generated by a two-dimensional map. It is clear that $MD_{p}$ is  non-smooth if, for some iteration, $x_{n+1}=x_n$ or $y_{n+1}=y_n$. Therefore, there are unbounded behaviours of the derivatives of $MD_{p}$ in those points independently of the dynamical behaviour of the system.

As mentioned in page 3531 in \cite{mancho CNS}, the heuristic law of  LDs is: ``singular contours of LDs correspond to invariant manifolds". We recall that in  a steady system with a saddle equilibrium, the stable (resp. unstable) manifold is defined as the set of points attracting to that equilibrium forwards (resp. backwards) in time \cite{palis}.
The concept of  {\it singular feature} is never presented by Mancho, Wiggins, and their co-workers in a precise way  which undoubtedly leads to ambiguity. Moreover, there is a lack of the definition of the structures it purports to detect.  For the reader's convenience, we consider the example discussed in page 3534 in \cite{mancho CNS}, specifically
\begin{equation}\label{Sl1}\left\{\begin{array}{lll}
x'=\lambda x\\
y'=-\lambda y\;\;\;\;{\rm with}\;\;\lambda>0.
\end{array}\right.\end{equation}
The structure of $M$ for $\tau=20$ and $t_{0}=0$ displays ``singular features" in the axes, the invariant manifolds of $(0,0)$ for system (\ref{Sl1}).
See FIG 1 (in this manuscript) or Figure 1 (c) or 2 (a) in \cite{mancho CNS}, (see also Section 2.1.2 in \cite{mancho CNS}).
\begin{figure}[ht]
\includegraphics[totalheight=2.8in]{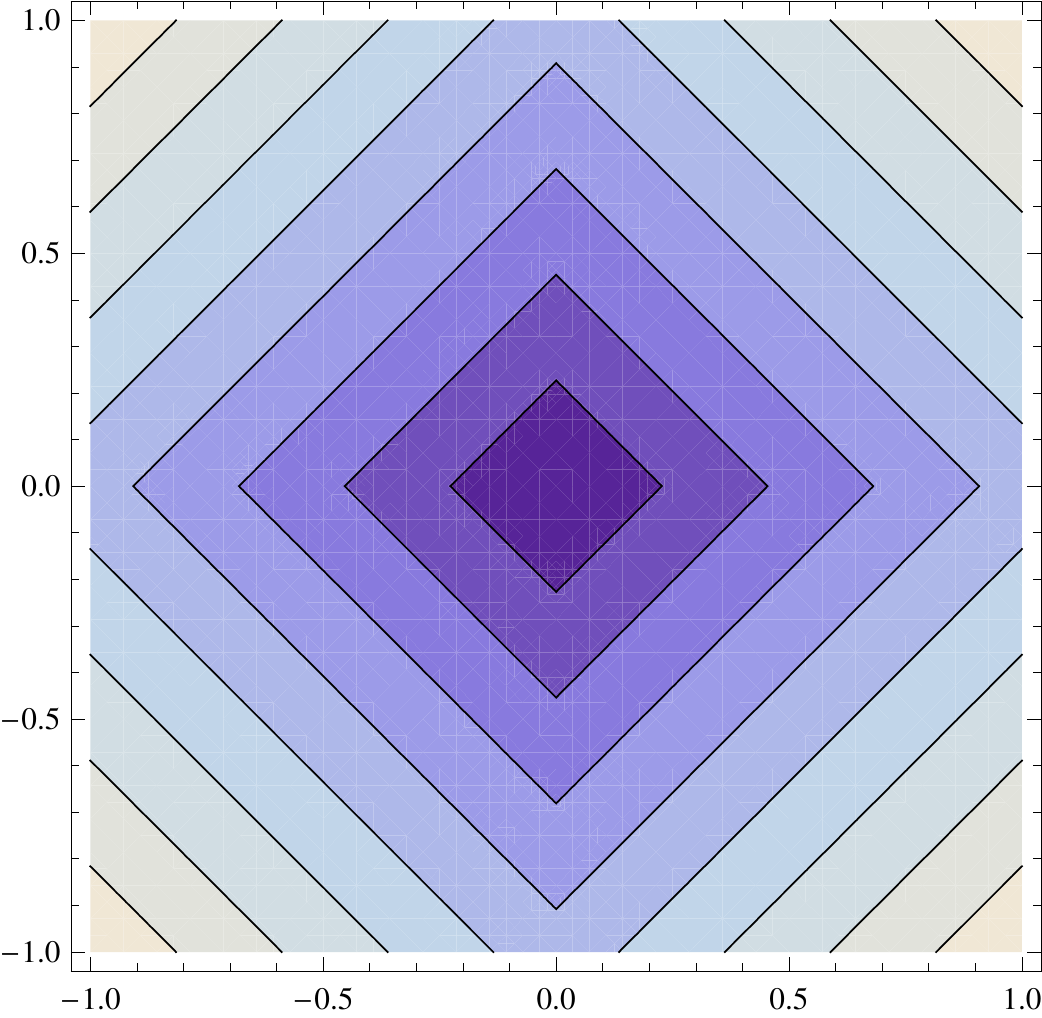}
\caption{ Contour-lines of $M$ in $(-1,1)\times(-1,1)$ associated with system (\ref{Sl1}) for $\lambda=1$ and $\tau=20$. The reader can check that this figure coincides with Fig. 1 (c) and Fig. 2 (a) in \cite{mancho CNS} (for $\tau=10$).}
\end{figure}

Recently, Lopesino {\it et al.} \cite{lopesino} have provided mathematical theorems to support the performance of the discrete version of LDs. Specifically, they stated the following result.
Let
\begin{equation}\label{S252}\left\{\begin{array}{lll}
x_{n+1}=\lambda x_{n}\\
y_{n+1}=\frac{1}{\lambda} y_{n}
\end{array}\right.\end{equation}
with $\lambda>1$.
\begin{theorem}\label{t10} (Theorem 1 in \cite{lopesino})
Consider a vertical line perpendicular to the unstable manifold of the origin. In particular consider an arbitrary point $x=\bar{x}$ and a line parallel to the $y$ axis passing through this point. Then the derivative of $MD_{p}$ with $p<1$, along this line becomes unbounded on the unstable manifold of the origin.
\end{theorem}
The proof  is just a straightforward check that $MD_{p}(x_{0},y_{0})=(|x_{0}|^{p}+|y_{0}|^{p})f(\lambda, N, x_{0},y_{0})$ with $f$ a smooth map.

We observe that Theorem \ref{t10}  is a consequence of the diagnostic itself because, as mentioned, $MD_{p}$ is non smooth if, for some iteration,
 \begin{equation}\label{contonta}
 x_{n+1}=x_n\;\;\; or\;\;\; y_{n+1}=y_n.
 \end{equation}
  Therefore, Theorem \ref{t10} gives no specific mechanism for detecting invariant manifolds: it only works for  systems satisfying (\ref{contonta}). We discuss this fact with a concrete example.

 Consider
\begin{equation}\label{S}\left\{\begin{array}{lll}
x_{n+1}= 2 x_{n}\\
y_{n+1}=\frac{1}{2} y_{n}+g(x_{n})
\end{array}\right.\end{equation}
with $g:\mathbb{R}\longrightarrow [0,\infty)$ a smooth function satisfying that $g(x)=0$ for all $x\not\in [0,1]$ and $g(x)>0$ if $x\in(0,1)$. The underlying map in (\ref{S}) is area preserving and $(0,0)$ is a global saddle point. Fix $\bar{x}>2$, it is straightforward to prove that given a line parallel to the $y$ axis passing through this point, the derivative of $MD_{p}$ with $p<1$, along this line becomes unbounded at $y=0$. Note that $(\bar{x},0)$ is not a point at the invariant manifolds of (\ref{S}). Thus, even under a slight modification of (\ref{S252}), Theorem \ref{t10} will give a false positive for an unstable or stable manifold. System  (\ref{S252}) also proves that Theorem \ref{t10} does not provide a mechanism to approximate invariant manifolds when $N\longrightarrow\infty$, (this property was mentioned in Section 2.1.2 in \cite{lopesino}).

Two remarks are in order:
\begin{itemize}
  \item In practical applications, Mancho, Wiggins, and their co-workers routinely use the contour-structure of the $M$-function (or $MD_{p}$ in discrete systems) to detect invariant manifolds, see \cite{mendoza-mancho PRL}-\cite{wiggins}. To  support the applicability of LDs, their arguments consist in proving that certain derivatives of $M$ or $MD_{p}$ are unbounded as $\tau\longrightarrow \infty$, (see the theorems in \cite{lopesino} or  Section 2.1.2 in \cite{mancho CNS}). Those arguments have no significance on the dynamical behaviour of the  systems or the contour structure of the $M$-function because $M$ and $MD_{p}$ are typically unbounded as $\tau\longrightarrow\infty$.

  \item In a recent work \cite{ruizherrera}, we gave precise families of examples in which LDs fail in the detection of  invariant manifolds. In our systems, the invariant manifolds were placed in the axes and the contour-lines of $M$ converge to horizontal lines as $\tau\longrightarrow \infty$, (the points inside these horizontal lines are indistinguishable for the  function $M$). In this scenario, one can just deduce  that the stable and unstable manifold are horizontal lines  independently of the definition of singular point under consideration.
\end{itemize}

 \section{ Lagrangian descriptors and Incompressible Flows}
 In this section we illustrate that the method of LDs fails in the  detection of invariant manifolds in 2D and 3D incompressible flows.\newline
 Consider
 \begin{equation}\label{Sl}\left\{\begin{array}{lll}
x'=f(x)\\
y'=- y f'(x)
\end{array}\right.\end{equation}
where $f:\mathbb{R}\longrightarrow \mathbb{R}$ is  of class $\mathcal{C}^{2}$ with $f(0)=0$ and satisfying the following conditions:
\begin{description}
  \item[C1] $f(x)>0$ if $x\in(0,\infty)$ and $f(x)<0$ if $x\in (-\infty,0)$.
  \item[C2] $f$ is bounded.
  \item[C3] $f'(0)\geq f'(x)>0$ for all $x\in\mathbb{R}$.
 \end{description}
We observe that (\ref{Sl}) has a global saddle at the origin where the $x$-axis (resp $y$-axis) is the un-stable (resp. stable) manifold. However, the structure of $M$ is smooth in a neighbourhood of the $y$-axis  as $\tau\rightarrow \infty$. Therefore, there are no singular features in a neighbourhood of the $y$-axis, (see FIG 2).

\begin{figure}[ht]
\includegraphics[totalheight=2.8in]{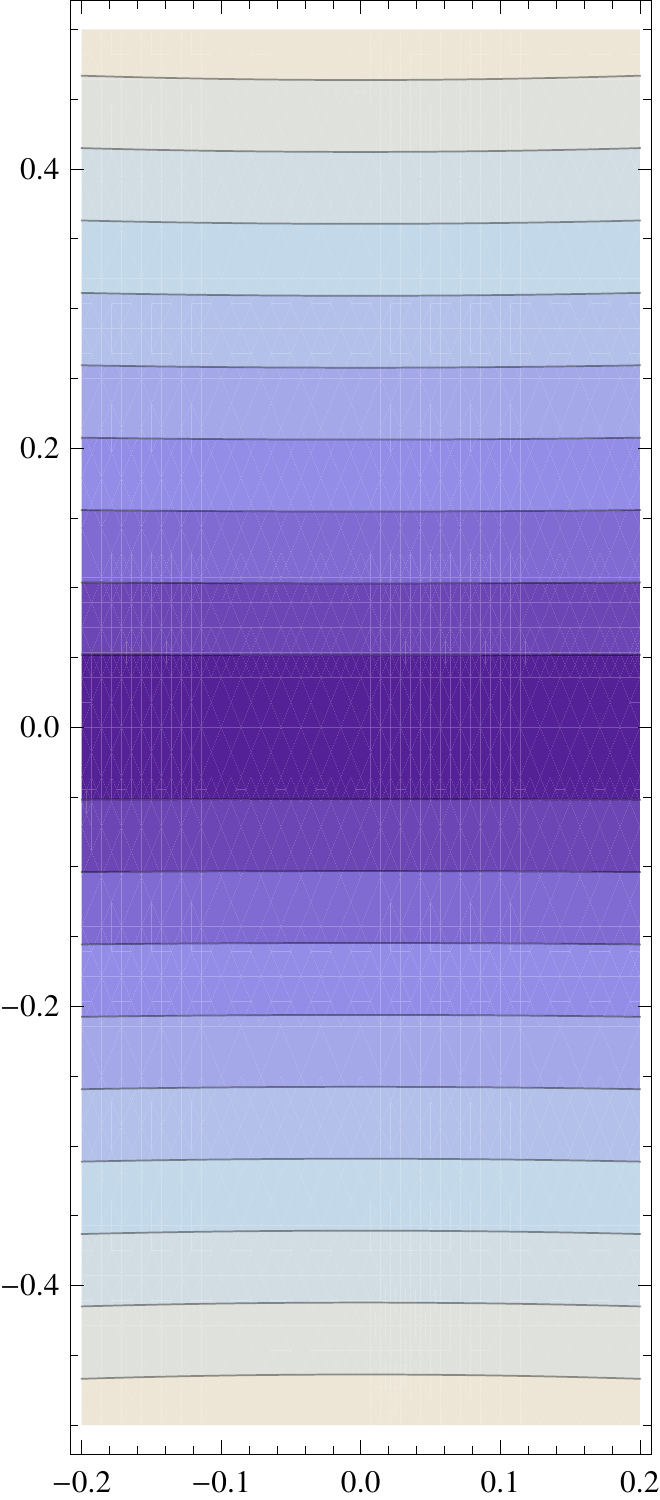}
\caption{Contour-lines of $M$ in $(-0.1,0.1)\times(-0.5,0.5)$ associated with system (\ref{Sl}) for $f(x)=\tanh x$ and $\tau=20$. One observes a smooth pattern in a neighbourhood of the $y$-axis .}
\end{figure}

Next, we give an analytic proof to justify this fact.
\begin{theorem}\label{tmanin} Assume that $f:\mathbb{R}\longrightarrow \mathbb{R}$ satisfies  {\bf C1-C3}
and
\begin{description}
  \item[C4] $f'(x)=k$ for all $x\in[-a,a]$ with $k,a>0$.
\end{description}
Then,  the contour lines of $M$ converge to horizontal lines as $\tau\longrightarrow \infty$ in a neighbourhood of the $y$-axis.
\end{theorem}
 \textbf{C4}  can be relaxed in a great deal, (e.g. we can consider $f''$ bounded in $\mathbb{R}$, $f''(x)>0$ in $(-\infty,0)$ and $f''(x)<0$ in $(0,\infty)$). However, the mathematical proof of Theorem \ref{tmanin} is much more involved. A prototypical function satisfying
\textbf{C1}-\textbf{C4} is
\begin{equation}\label{f}f(x)= \left\{ \begin{array}{lcc}
             \arctan(x+1)-1 &  if  & x \leq -1 \\
             \\ x &  if & -1 < x < 1 \\
             \\ \arctan(x-1)+1 &  if  & x \geq 1.
             \end{array}
   \right.\end{equation}

To prove Theorem \ref{tmanin}, we need the following lemma.

\begin{lemma}\label{L2} Assume the conditions  of Theorem \ref{tmanin}.  Given $(x_{0},y_{0})\in[-a,a]\times (0,\infty)$ then
$$\lim_{\tau\longrightarrow \infty}\frac{M(\tau;x_{0},y_{0})}{M(\tau;0,y_{0})}=1.$$
\end{lemma}
\textbf{Proof.}
First of all, we note that the solution of (\ref{Sl}) with initial condition $(x_{0},y_{0})$ is given by
$$(x(t;x_{0},y_{0}),y(t;x_{0},y_{0}))=(x(t;x_{0}),y_{0}e^{-\int _{0}^{t}f'(x(s;x_{0}))ds})$$
where $x(t;x_{0})$ is the solution of $x'=f(x)$ with initial condition $x_{0}$.
By simple computations and using  \textbf{C2} and \textbf{C3},
$$M(\tau;0,y_{0})=y_{0} (e^{f'(0)\tau}-e^{-f'(0)\tau}),$$
$$M(\tau;x_{0},y_{0})=\int _{-\tau}^{\tau}\sqrt{x'(t;x_{0})^{2}+y_{0}^2 e^{-2\int_{0}^{t}f'(x(s;x_{0}))ds}f'(x(t;x_{0}))^{2}}dt $$
$$\leq \int _{0}^{\tau}\sqrt{f(x(t;x_{0}))^{2}+y_{0}^2 f'(0)^2}dt+ \int _{-\tau}^{0}\sqrt{f(x(t;x_{0}))^{2}+y_{0}^2 e^{-2f'(0)t}f'(0)^{2}}dt $$
$$\leq \int _{0}^{\tau}\sqrt{m^{2}+y_{0}^2 f'(0)^2}dt+ \int _{-\tau}^{0}\sqrt{m^{2}+y_{0}^2 e^{-2f'(0)t}f'(0)^{2}}dt$$

for all $(x_{0},y_{0})$ with $y_{0}>0$ where $m$ is a bound of $f$.
Therefore,

$$\limsup_{\tau\longrightarrow \infty}\frac{M(\tau,x_{0},y_{0})}{M(\tau,0,y_{0})}\leq 1$$
because of $$ \lim _{\tau\longrightarrow \infty}\frac{\tau \sqrt{m^{2}+y_{0}^2 f'(0)^2}}{y_{0} (e^{f'(0)\tau}-e^{-f'(0)\tau})}=0,$$$$ \lim _{\tau\longrightarrow \infty}\frac{\int _{-\tau}^{0}\sqrt{m^{2}+y_{0}^2 e^{-2f'(0)t}f'(0)^{2}}dt}{y_{0}(e^{f'(0)\tau}-e^{-f'(0)\tau})}=1.$$
The second limit is a simple application of the L'H\^{o}pital rule.

By \textbf{C1} and \textbf{C4}, $$-a\leq x(t;x_{0})\leq a$$ for all $t<0$ and $-a\leq x_{0}\leq a$. Thus, $f'(x(t;x_{0}))=f'(0)=k$ for all $t<0$ and $-a\leq x_{0}\leq a$.
Now, we conclude that
$$\liminf_{\tau\longrightarrow \infty}\frac{M(\tau;x_{0},y_{0})}{M(\tau;0,y_{0})}\geq 1$$
since
$$M(\tau;x_{0},y_{0})=\int _{-\tau}^{\tau}\sqrt{x'(t;x_{0})^{2}+y_{0}^2 e^{-2\int_{0}^{t}f'(x(s;x_{0}))ds}f'(x(t;x_{0}))^{2}}dt $$
$$\geq \int _{-\tau}^{0}\sqrt{x'(t;x_{0})^{2}+y_{0}^2 e^{-2\int_{0}^{t}f'(x(s;x_{0}))ds}f'(x(t;x_{0}))^{2}}dt\geq y_{0}( e^{ f'(0) \tau}-1).$$ \qed

\textbf{Proof of Theorem \ref{tmanin}.} For each $x_{0}\in[-a,a]$, we define $y_{0}^{\tau}$ such that $(x_{0},y_{0}^{\tau})\in \gamma ^{\tau}_{(0,y_{0})}$,  where $\gamma_{(0,y_{0})}^{\tau}$ is the contour-line of $M$ passing through $(0,y_{0})$ for $\tau$.
First we observe that for all $\tau>0$ large enough, $y_{0}^{\tau}$ exists since
$$M(\tau;x_{0},0)\leq 2 m \tau < M(\tau;0,y_{0})=y_{0}(e^{f'(0)\tau}-e^{-f'(0)\tau})$$
where $m$ is a bound of $f$ and
$$ M(\tau;x_{0},2 y_{0})\geq 2 y_{0}(e^{f'(0)\tau}-1)\geq M(\tau;0,y_{0}). $$

 By Lemma \ref{L2},
$$M(\tau,x_{0},y_{0})=M(\tau,0,y_{0})+f(\tau,x_{0},y_{0})$$ with
$$\lim_{\tau\longrightarrow \infty} \frac{f(\tau,x_{0},y_{0})}{e^{f'(0)\tau}-e^{-f'(0)\tau}}=0$$
for all $(x_{0},y_{0})\in K$. Therefore, by using that $M(\tau,0,y_{0})=M(\tau,x_{0},y_{0}^{\tau})$, we have that

$$y_{0}(e^{f'(0)\tau}-e^{-f'(0)\tau})= y_{0}^{\tau}(e^{f'(0)\tau}-e^{-f'(0)\tau})+f(\tau,x_{0},y_{0}^{\tau}).$$
From this expression,
$$y_{0}=y_{0}^{\tau} + \frac{f(\tau,x_{0},y_{0}^{\tau})}{e^{f'(0)\tau}-e^{-f'(0)\tau}}$$
and clearly, $\lim_{\tau\longrightarrow \infty}y_{0}^{\tau}=y_{0}$.

Collecting all the information, we have proved that

\begin{equation}\label{condition1}
\gamma_{(0,y_{0})}^{\tau}\cap ([-a,a]\times (0,\infty))\rightarrow [-a,a]\times \{y_{0}\}.
\end{equation}

In other words, the contour lines in a neighbourhood of the stable manifold tend to horizontal segments. \qed\newline

Next we show the inability of the method of LDs to detect invariant manifolds  in linear systems.
\begin{theorem}\label{t2}
Consider
\begin{equation}\label{S25}\left\{\begin{array}{lll}
x_{1}'=\lambda_{1} x_{1}\\
x_{2}'=\lambda_{2} x_{2}\\
\vdots\\
x'_{N-1}=\lambda_{N-1}x_{N-1}\\
x'_{N}=-\lambda_{N} x_{N}
\end{array}\right.\end{equation}
with $\lambda_{i}>0$ and $\lambda_{N}>\max\{\lambda_{i}:i=1,...,N-1\}$. Then
the contour-surfaces of $M$ in a neighbourhood of the $x_{N}$-axis converge to hyper-planes parallel to $x_{N}=0$ as $\tau\longrightarrow \infty$.
\end{theorem}

\textbf{Proof.} Denote by $\{e_{i}:i=1,...,N\}$ the usual basis of $\mathbb{R}^{N}$. It is clear that
$$M(\tau;x_{i}e_{i})=x_{i}(e^{\lambda_{i}\tau}-e^{-\lambda_{i}\tau})$$ and
$$M(\tau;x_{N}e_{N})\leq M(\tau;(x_{1},...,x_{N}))\leq  M(\tau;x_{1}e_{1})+\ldots + M(\tau;x_{N}e_{N})$$
for all $(x_{1},\ldots,x_{N})\in\mathbb{R}^{N}$. Therefore, for $(x_{0},...,x_{N})\in \mathbb{R}^{N-1}\times (0,\infty)$, we have that
$$\lim_{\tau\longrightarrow \infty}\frac{M(\tau;(x_{1},...,x_{N}))}{M(\tau;x_{N}e_{N})}=1$$ for all $(x_{1},...,x_{N})\in K $. By repeating the proof of the previous theorem, we get that
\begin{equation}\label{condition1}
\gamma_{x_{N}e_{N}}^{\tau}\cap ([-a,a]^{N-1}\times (0,\infty))\rightarrow [-a,a]^{N-1}\times \{y_{N}\}
\end{equation}
with $a>0$. We have denoted by $\gamma_{x_{N}e_{N}}^{\tau}$  the contour-surface of $M$ for $\tau$ passing through $x_{N}e_{N}$. That is, the contour surfaces of $M$ in a neighbourhood of the stable manifolds tend to a hyperplane parallel to $x_{N}=0$.\qed\newline
 We observe that Theorem \ref{t2} can be applied  in incompressible flows, for instance $\lambda_{1}=\lambda_{2}=1$ and $\lambda_{3}=2$.

\section{Minima of $M$ do not detect invariant manifolds}
Craven and Hern\'{a}ndez in \cite{prl}  have employed a modified version of the method of Lagrangian Descriptors  to study the dynamical skeleton of
Langevin equation
\begin{equation}\label{1}
m \ddot{q}=-\gamma \dot{q} -\frac{\partial V}{\partial q}(q,t)+ \sqrt{2\sigma}\xi_{a}(t)\;\;\; \;\;\;\;m, \sigma, \gamma > 0
\end{equation}
where  $\xi_{a}(t)$ is a noise input and $V(q,t)$ represents the underlying potential.  The authors claim that the stable  manifold (at $t_{0}$) in (\ref{1}) is determined {\it holding the coordinate $q=C$ constant} (transversal direction) {\it and minimizing with respect to $\dot{q}_{0}$}
\begin{equation}\label{Lf}
L_{f}((C,\dot{q}_{0}),t_{0})_{\tau}=\int _{t_0}^{t_{0}+\tau}\|{\bf{\dot{q}}_{c}}(C,\dot{q}_{0},t_{0},t)\|dt
\end{equation}
\noindent where ${\bf{q_{c}}}(C,\dot{q}_{0},t_{0},t)=(q_{c}((C,\dot{q}_{0}),t_{0},t),\dot{q}_{c}((C,\dot{q}_{0}),t_{0},t))$ is the solution of (\ref{1}) satisfying ${\bf q_{c}}(C,\dot{q}_{0},t_{0},t_{0})=(C,\dot{q}_{0})$. See Figure 2 (g) in \cite{prl} for a pictorial explanation of this heuristical law. In that figure, one observes that the authors associate the global minima of $L_{f}$ with the stable manifold, see also \cite{craven and hernandez, craven and hernandez 1}. The objective of this section is to show that this recipe does not provide a suitable mechanism to detect  invariant manifolds in Langevin equations.

First we consider
 \begin{equation}\label{2}
 \left\{\begin{array}{lll}
x'=-x+y\\
y'=h(y)
\end{array}\right.
\end{equation}
with
$$h(y)=
\left\{\begin{array}{lll}
y-y^2&if\;\;y\geq0\\
y+y^2&if\;\;y<0.
\end{array}\right.
$$
In (\ref{2}), the $x$-axis separates the basis of attraction of  equilibria $(1,1)$ and $(-1,-1)$. Holding $x=x_{0}$ constant with $x_{0}>1$,
the global minimum of $L_{f}((x_{0},y),0)_{\tau}$ is  attained   at (aprox.) $1$, not $0$ (the location of the stable manifold), (see FIG. 3 left). In fact $\frac{d}{dy}L_{f}((x_{0},y),0)_{\tau}\Big|_{y=0}=\frac{-e^{\tau}+e^{-\tau}}{2}.$ Notice that $(x_{0},1)$ is not a point at the  invariant manifolds. This is the typical dynamical situation studied in \cite{prl}, a manifold which separates reactive and non-reactive trajectories. We have taken system (\ref{2}) because the solutions and all quantities are computable with simple mathematical skills. A particular counter-example of type (\ref{1}) with similar pathologies is
\begin{equation}\label{r1}
\ddot{q}=-\dot{q}+0.1 q(1-q^{2}).
\end{equation}

This equation has two attractors, namely $(1,0), (-1,0)$, and the origin is a saddle point. The graph of $L_{f}((1.1,\dot{q}),0)_{\tau}$ is illustrated in FIG.3 right. The global minimum $L_{f}((1.1,\dot{q}),0)_{\tau}$ is attained in the interval  $(-0.1,0.1)$ for all $\tau$. However, $(1.1,6.17)$ (aprox.) is the point of the stable manifold.
\begin{figure}[ht]
\includegraphics[totalheight=1.8in]{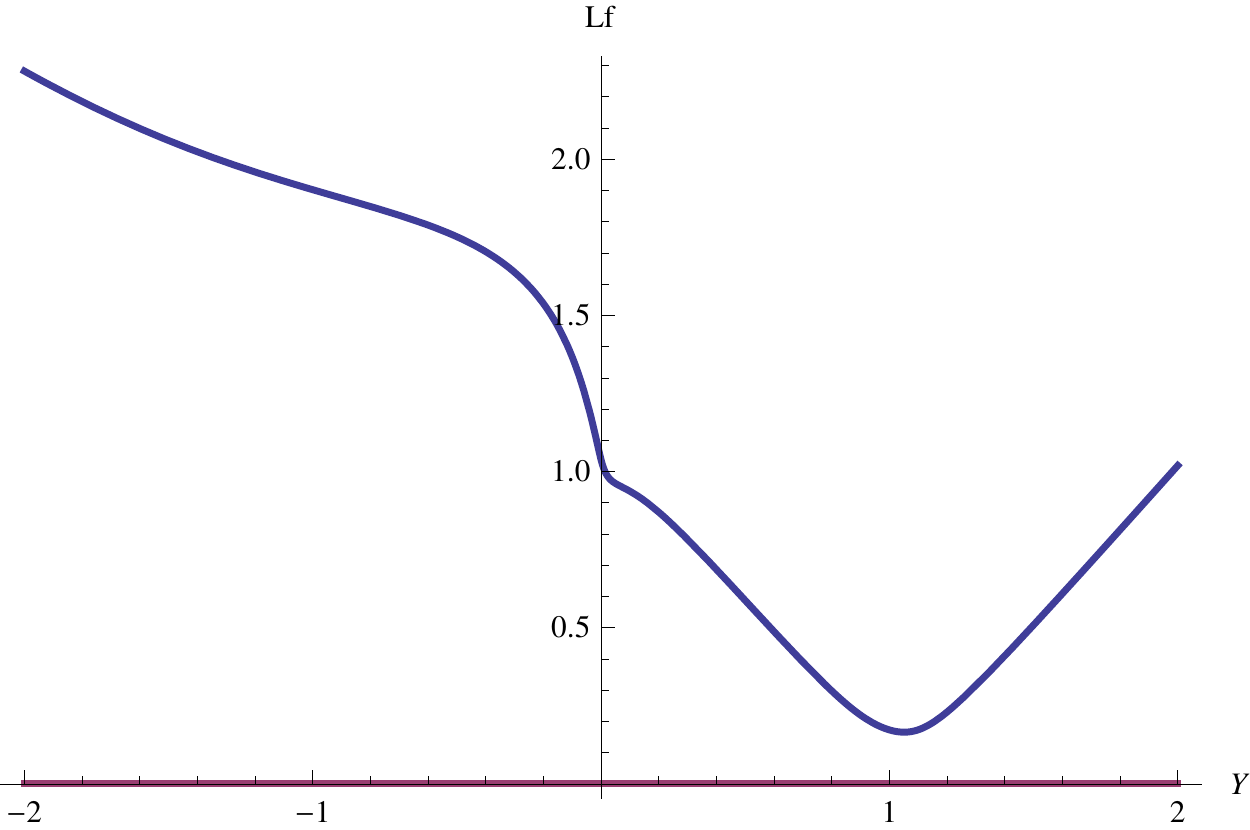}
\includegraphics[totalheight=1.8in]{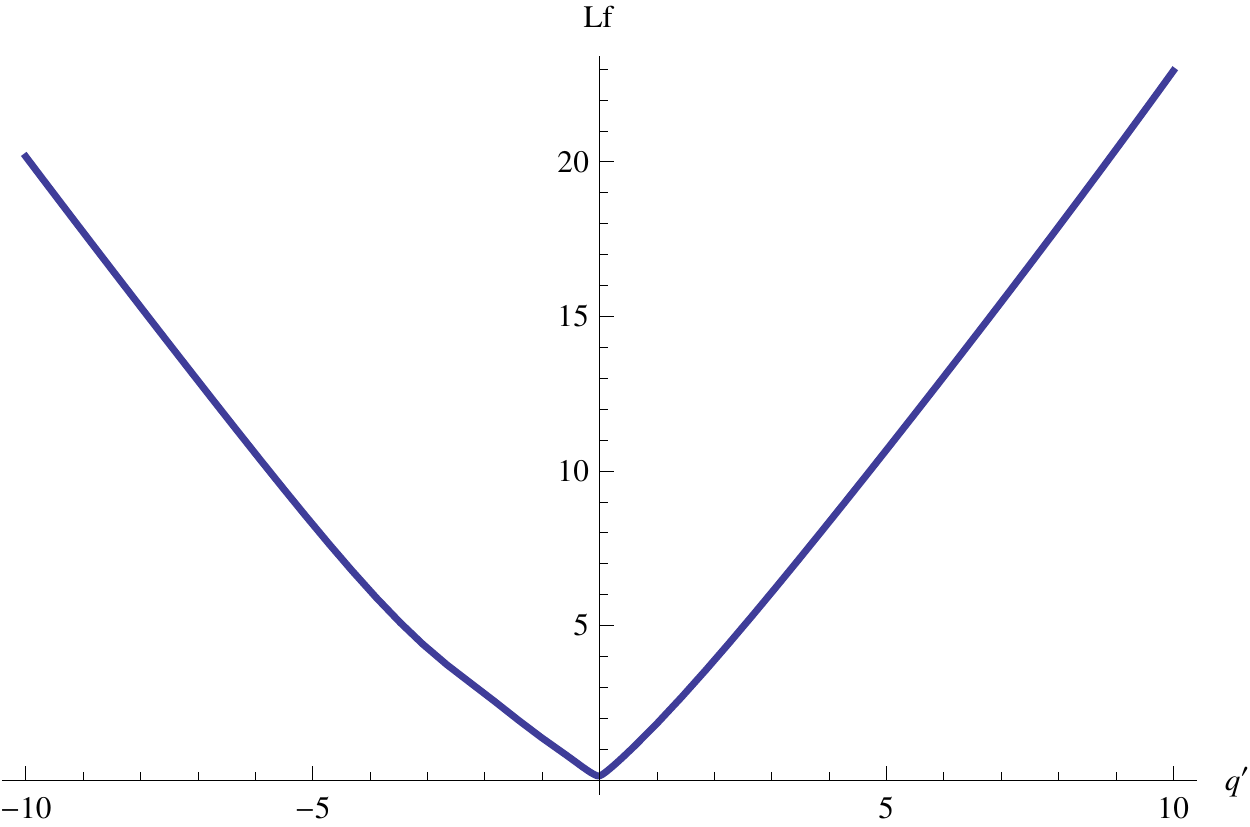}
\caption{Representation of $L_{f}((1.1,y),0)_{2}$ in $y\in[-2,2]$ for system (\ref{2}) (left) and $L_{f}((1.1,q'),0)_{20}$ in $[-10,10]$ for equation (\ref{r1}).
   The global minimum is not attained at the stable manifold.}
\end{figure}

\section{Discussion}

The objective of this paper has been to illustrate that  singularities of the contour-lines of $M$ do not generally detect invariant manifolds in incompressible flows. We have constructed systems in which the stable and unstable manifolds are placed in the axes but the contour lines of $M$ converge to horizontal lines as $\tau\longrightarrow\infty$ (the points inside the segments are indistinguishable for the function $M$). Hence, if the contour-lines of $M$ for large values were useful in the detection of  invariant manifolds, one would deduce that the stable and unstable manifolds are horizontal lines (independently of the definition of singularity under consideration).

In real applications,  LDs involve finite (large) values of $\tau$. As mentioned above, from a mathematical point of view, the contour-structure of the $M$ function is smooth for all $\tau$. Therefore, the possible ``singularities" would appear at $\tau=\infty$. Notice that if the contour-lines of $M$ converge to horizontal lines as $\tau\longrightarrow\infty$, horizontal lines  will appear (in our computer) for $\tau$ large enough.

 Except for the simple system
 \begin{equation}\label{2}
 \left\{\begin{array}{lll}
x'=-\lambda x\\
y'=\lambda y,
\end{array}\right.
\end{equation}
there are no  mathematical basis to justify the applicability  of  LDs. Therefore, the conclusions inferred    in \cite{mendoza-mancho PRL}-\cite{wiggins} need some type of updating.  Of course, the $M$-function  always creates certain patterns when plotted over initial conditions. However, as emphasized in this manuscript, this fact does not imply that the output has dynamical significance.

 In \cite{balibrea response}, Mancho, Wiggins, and their co-workers claim that criticisms along the lines of \cite{ruizherrera} or the present paper are not relevant because the conclusions of LDs are not based on the contour-structure of the $M$-function. This claim is inaccurate. It seems that they misrepresent what they have actually done. To check that Mancho, Wiggins and their co-workers routinely use the contour-structure of $M$ in their work, the reader can see
 figures 1-4  in \cite{mendoza-mancho PRL}; figures 1-3, 5,7,9 in \cite{camara 1} (the contour-lines of the derivative of $M$ are missing); see figures of section 4 (applications) in \cite{lopesino}; see figures 1-5 and section 4.1 in \cite{mancho avion}  (this is their last application) and so on.
 In \cite{balibrea response}, they suggested  that one has to look for points  at which certain directional derivatives of $M(x_0, t_0, \tau )$
do not exist, inspired by Lopesino {\it et al} \cite{lopesino}. As mentioned in Section II, the results in \cite{lopesino}  give no a recipe to detect invariant manifolds for general systems. Moreover, those results were given to extract information on the contour-lines of $MD_{p}$ (see the applications and introduction in \cite{lopesino}). The function $M(x_0, t_0, \tau )$ is smooth (except at equilibria) and is usually unbounded as $\tau\longrightarrow\infty$.  The absence of discontinuity of $M$, as in the systems discussed in \cite{balibrea response}, appears in most  systems everywhere. For instance, in system
 \begin{equation}\label{3}
 \left\{\begin{array}{lll}
x'=- y\\
y'= x,
\end{array}\right.
\end{equation}
the $M$-function is given by
$M(\tau,x,y)=2\tau \sqrt{x^2+y^2}$ and the partial derivatives are unbounded if $x\not=0$ or $y\not=0$ when $\tau\longrightarrow\infty$.

\section*{Funding Statement}
 This research was supported by Spanish grant MTM2014-56953-P.


\begin{thebibliography}{777}
\bibitem{mendoza-mancho PRL}{\sc
C. Mendoza and A. M. Mancho}, {\it
 The hidden geometry
of ocean flows}, Physical Review Letters 105 (2010) 038501.
\bibitem{mancho CNS}{\sc
 A. M. Mancho,
S. Wiggins, J. Curbelo, and C. Mendoza,} {\it Lagrangian Descriptors:
A Method for Revealing Phase Space Structures of General Time
Dependent Dynamical Systems}, Communications in Nonlinear Science
and Numerical Simulation 18 (2013), 3530--3557.
\bibitem{lopesino}{\sc C. Lopesino, F. Balibrea, S. Wiggins,  A.M. Mancho}, {\it Lagrangian Descriptors for Two Dimensional, Area Preserving Autonomous and Nonautonomous Maps}, Communications in Nonlinear Science and Numerical Simulation 27  (2015), 40--51.

\bibitem{mancho libro}{\sc A. M. Mancho, J. Curbelo, S. Wiggins, V.J. Garcia-Garrido, and C. Mendoza},   Beautiful Geometries Underlying Ocean Nonlinear Processes Chapter in the book.   A Voyage Through Scales. Eds. Günter Blöschl, Hans Thybo, Hubert Savenije, Lois Lammerhuber. Publishers European Geophysical Union and Edition Lammerhuber (2015).

\bibitem{smith}{\sc M.L. Smith and A.J. McDonalds}, {\it
A quantitative measure of polar vortex strength
using the function M},    J. Geophys Res: Atmos 119 (2014), 5966--5986.
\bibitem{rempel}{\sc E.L. Rempel, A. Chian, A. Brandenburg, P. Munoz, and S. Shadden}, {\it Coherent structures and the saturation of a nonlinear dynamo}, J. Fluid Mech. 729 (2013), 309--329.
\bibitem{mendoza npg2014}{\sc C. Mendoza, A. M. Mancho, and S. Wiggins,} {\it Lagrangian Descriptors and
the Assesment of the Predictive Capacity of Oceanic Data Sets},
Nonlinear Processes in Geophysics 21 (2014), 485--501.

\bibitem{mancho npg2014 nekhoroshev}{\sc S. Wiggins and A. M Mancho,} {\it Barriers to tranport in
aperiodically time-dependent two dimensional velocity fields:
Nekhoroshev's Theorem and 'Nearly Invariant' Tori},  Nonlinear
Processes in Geophysics 21 (2014), 165--185.
\bibitem{camara 1}{\sc
 A. de la C{\'a}mara, R. Mechoso, A. M. Mancho, E. Serrano, and K.
Ide, } {\it  Quasi-horizontal transport within the Antarctic polar
night vortex: Rossby wave breaking evidence and Lagrangian
structures}, Journal of the Atmospheric Sciences 70 (2013)
2982--3001.

\bibitem{mendoza npg2012}{\sc
C. Mendoza and A. M. Mancho, } {\it The Lagrangian description of
aperiodic flows: a case study of the Kuroshio Current}, Nonlinear
Processes in Geophysics 19  (2012), 449--472.
\bibitem{camara 2}{\sc
A. de la C{\'a}mara, A. M. Mancho, K. Ide, E. Serrano, and C.R. Mechoso,}
{\it Routes of transport across the Antarctic polar vortex in the
southern spring}, Journal of the Atmospheric Sciences 69  (2012)
753--767.


\bibitem{mendoza npg 2010}{\sc
C. Mendoza, A. M. Mancho, and M.-H. Rio},{\it The turnstile mechanism
across the Kuroshio current: analysis of dynamics in altimeter
velocity fields}, Nonlinear Proc. Geoph 17 (2010),  103-111.


\bibitem{mancho avion}{\sc V. J. Garc\'{i}a-Garrido, A.M. Mancho, S. Wiggins, and C. Mendoza}, {\it A dynamical systems approach to the surface search for debris associated with the disappearance of flight MH370}, Nonlinear Processes in Geophysics, 22 (2015), 701.

\bibitem{Guha}{\sc A. Guha, C.R. Mechoso, C.S. Konor, and R.P. Heikes,} {\it  Modeling Rossby Wave Breaking in the Southern Spring Stratosphere}, Journal of the Atmospheric Sciences 73 (2016), 393--406.
\bibitem{nueva tonteria}{\sc R. Vortmeyer-Kley, U. Grave, and U. Feudel,} {\it Detecting and tracking eddies in oceanic flow fields: A Lagrangian Descriptor based on the modulus of vorticity},  Nonlinear Processes in Geophysics, 23 (2016), 159.

\bibitem{wiggins}{\sc V.J. Garcia-Garrido, A. Ramos, A.M. Mancho, J. Coca, S. Wiggins}, {\it A dynamical system perspective for a real-time response to a marine oil spill}, Marine Population Bulletin (in press)

\bibitem{mancho chaos}{\sc
J. A. Madrid and A. M. Mancho,}
{\it Distingueshed trajectories in time dependent vector fields}, Chaos 19 (2009), 013111.
\bibitem{ruizherrera}{\sc A. Ruiz-Herrera}, {\it Some examples related to the method of Lagrangian descriptors}, Chaos  25 (2015), 063112.
\bibitem{palis} {\sc J.J. Palis and W. De Melo},  Geometric theory of dynamical systems: an introduction. Springer Science and Business Media (2012).
\bibitem{haller M}{\sc G. Haller},  {\it Non-objectivity of the M function and other thoughts}, http://www.nonlin-processes-geophys-discuss.net/npg-2016-16/
\bibitem{chaos intro} {\sc T. Peacock, G. Froyland, and G. Haller,} {\it  Introduction to Focus Issue: Objective Detection of Coherent Structures.} Chaos 25 (2015), 7201.

\bibitem{craven and hernandez} {\sc G.T. Craven and R. Hernandez,} {\it  Deconstructing field-induced ketene isomerization through Lagrangian descriptors}, Physical Chemistry Chemical Physics (2016)


\bibitem{prl} {\sc G. Craven and R. Hernandez}, {\it Lagrangian Descriptors of Thermalized Transition States on Time-Varying Energy Surfaces}, Physical Review Letters 115 (2015),   148301.
\bibitem{craven and hernandez 1}{\sc A. Junginger and R. Hernandez,} {\it  Uncovering the Geometry of Barrierless Reactions Using Lagrangian Descriptors} The Journal of Physical Chemistry B 120 (2016), 1720–1725.
\bibitem{balibrea response}{\sc F. Balibrea-Iniesta, J. Curbelo, V.J. Garcia-Garrido, C. Lopesino, A.M. Mancho, C. Mendoza, and S. Wiggins,} Response to:" Limitations of the Method of Lagrangian Descriptors"[arXiv: 1510.04838]. arXiv preprint arXiv:1602.04243. (2016)
 \bibitem{haller}{\sc G. Haller}, {\it Lagrangian coherent structures},  Annual Review of Fluid Mechanics 47 (2015), 137--162.
\bibitem{allshouse}{\sc M.R. Allshouse and T. Peacock}, {\it Refining finite-type Lyapunov exponents ridges and the challengues of classifying them}, Chaos 25 (2015), 087410.
 \bibitem{Froyland1}{\sc G. Froyland}, {\it An analytical framework for identifying finite-time coherent sets in time dependent dynamical systems}, Physica D 250 (2013), 1--19.
 \bibitem{Haller2}{\sc D. Karrasch, F. Huhn, and G. Haller}, {\it Automated detection of coherent Lagrangian vortices in two dimensional unsteady flows}, Proceedings of the Royal Society of London A 471 (2014), 20140639.
  \bibitem{Froyland2} {\sc G. Froyland}, {\it Dynamic isoperimetry and the geometry of Lagrangian Coherent Structures}, Nonlinearity 25 (2015), 087406.
  \bibitem{peacock} {\sc T. Peacock and G. Haller}, {\it Lagrangian Coherent Structures: The hidden skeleton of fluid flows}, Physics Today 66 (2013), 41.



\end{thebibliography}
\end{document}